\documentclass[12pt,a4paper]{article}
\usepackage{amssymb, amscd, amsthm, amsmath,latexsym, amstext,mathrsfs}
\usepackage[all]{xy}
\usepackage[dvips]{graphicx}

\newtheorem{thm}{Theorem}

\newtheorem{lemma}{Lemma}
\newtheorem{prop}{Proposition}

\theoremstyle{remark}
\newtheorem{remark}{\textbf{Remark}}

\theoremstyle{definition}
\newtheorem{defn}{Definition}

\newtheorem{examples}{Example}
\newtheorem{prob}{Problem}

\newcommand{\empha}[1]{\textbf{\emph{#1}}}

\def\deg{\mathrm{deg}}

\begin{document}
\title{Weak Approximation over Function Fields of Curves over Large or Finite Fields}
\author{Yong HU\footnote{Math\'{e}matiques, B\^{a}timent 425, Universit\'{e} Paris-Sud, 91405,
Orsay Cedex,  France,\ \ \ e-mail: yong.hu2@math.u-psud.fr}}
\date{}

\maketitle

\begin{abstract}
Let $K=k(C)$ be the function field of a curve over a field $k$ and
let $X$ be a smooth, projective, separably rationally connected
$K$-variety with $X(K)\neq\emptyset$. Under the assumption that $X$
admits a smooth projective model $\pi: \mathcal{X}\to C$, we prove
the following weak approximation results: (1) if $k$ is a large
field, then $X(K)$ is Zariski dense; (2) if $k$ is an infinite
algebraic extension of a finite field, then $X$ satisfies weak
approximation at places of good reduction; (3) if $k$ is a
nonarchimedean local field and $R$-equivalence is trivial on one of
the fibers $\mathcal{X}_p$ over points of good reduction, then there
is a Zariski dense subset $W\subseteq C(k)$ such that $X$ satisfies
weak approximation at places in $W$. As applications of the methods,
we also obtain the following results over a finite field $k$: (4) if
$|k|>10$, then for a smooth cubic hypersurface $X/K$, the
specialization map $X(K)\longrightarrow \prod_{p\in
P}\mathcal{X}_p(\kappa(p))$ at finitely many points of good
reduction is surjective; (5) if $\mathrm{char}\,k\neq 2,\,3$ and
$|k|>47$, then a smooth cubic surface $X$ over $K$ satisfies weak
approximation at any given place of good reduction.
\end{abstract}

\tableofcontents

\

\noindent {\bf Convention.} In this paper, a \empha{variety} means a
separated geometrically integral scheme of finite type over a field.
Somewhat inconsistently, a \empha{curve} means a 1-dimensional
separated geometrically connected and geometrically reduced scheme
of finite type over a field. The words ``\empha{vector bundle}''
will be used interchangeably as ``locally free sheaf of finite
constant rank''.

\section{Introduction}

Given a variety over a topological field $F$, the topology of $F$
induces a natural topology on the set of rational points. We call
this topology the \empha{$F$-topology}. If $K\subseteq F$ is a
subfield and $X$ is a $K$-variety, the set $X(F)$ of $F$-points
carries an $F$-topology via the natural identification
$X(F)=X_F(F)$, where $X_F=X\times_KF$ is the $F$-variety obtained by
base extension of the $K$-variety $X$. Let $K$ be a number field or
a function field of a curve. The completion $K_v$ at each place $v$
of $K$ is a topological field containing $K$. The $K_v$-topologies
will be also called \empha{$v$-adic topologies}. Given a set
$\Omega$ of places of $K$, we say $X$ satisfies \empha{weak
approximation} at places in $\Omega$ if the image of the diagonal
map $X(K)\to\prod_{v\in \Omega}X(K_v)$ is dense, where $\prod_{v\in
\Omega}X(K_v)$ is given the product topology of the $v$-adic
topologies. The classical approximation theorem claims that the
projective line $\mathbb{P}^1$ satisfies weak approximation at all
places. In higher dimension and
 over the function field of a curve defined over
an algebraically closed field, it is expected that weak
approximation holds for \empha{separably rationally connected}
varieties. Roughly speaking, these are varieties containing plenty
of geometric rational curves that may be deformed ``very freely''.
(See Definition$\;$\ref{defn2p1} for the precise definition.)

In this paper, we are interested in the following problem.

\begin{prob}\label{prob1p1}
 Let $k$ be a field, $C$ a smooth projective (geometrically connected) curve over
 $k$, $K=k(C)$ the function field of $C$, and $X$ a smooth,
 projective, separably rationally connected $K$-variety. Assume that
 $X$ has a smooth projective model, i.e., there is
 a surjective morphism $\pi: \mathcal{X}\to C$ from a smooth projective $k$-variety $\mathcal{X}$ to $C$
 with generic fiber $X/K$. In this
context, rational points of $X$ correspond to sections of $\pi:
\mathcal{X}\to C$ and weak approximation amounts to the existence of
sections passing through prescribed points in the fibers with
prescribed jets.

Consider the following questions:

(EX) (Existence of sections) Is $X(K)$ nonempty? Or equivalently,
does $\pi$ admit a section $\sigma: C\to \mathcal{X}$?

(SP) (Surjectivity of specialization map) Suppose that $\pi$ has a
section $\sigma: C\to\mathcal{X}$. Let $P=\{\,p_1,\dotsc, p_n\,\}$
be a finite set of closed points of $C$ such that each fiber
$\mathcal{X}_i=\mathcal{X}_{p_i}$ is smooth and separably rationally
connected (a point or a place with this property will be called
\empha{of good reduction}). Given rational points $r_i\in
\mathcal{X}_i(\kappa(p_i))$ in the fibers, is there a section $s:
C\to \mathcal{X}$ of $\pi$ such that $s(p_i)=r_i$ for $i=1,\dotsc,
n$?

(WA) (Weak approximation) Suppose that $X(K)\neq \emptyset$. Does
the generic fiber $X$ satisfy weak approximation at a certain
collection of places?
\end{prob}

When $k$ is algebraically closed, a series of results have been
established. Firstly, Koll\'{a}r, Miyaoka and Mori \cite{KMM} proved
that question (SP) has a positive answer in this case. Later,
question (EX) was solved by Graber, Harris and Starr \cite{GHS} in
characteristic 0 and soon generalized to positive characteristic by
de Jong and Starr \cite{deJS}. As for question (WA),
Colliot-Th\'{e}l\`{e}ne and Gille \cite{CTG} proved that if $X$
belongs to one of some special classes of varieties, including conic
bundles and del Pezzo surfaces of degree $\ge 4$, then $X$ satisfies
weak approximation at all places. For arbitrary separably rationally
connected varieties, Hassett and Tschinkel \cite{HaTs1} proved weak
approximation at places of good reduction, generalizing the case of
cubic surfaces proved earlier by Madore \cite{Mad}. In the same
paper, they conjectured that weak approximation at places of bad
reduction also holds. This conjecture is open even for cubic
surfaces. Hassett and Tschinkel \cite{HaTs2} confirmed weak
approximation for cubic surfaces with not too bad reductions. A
similar result was obtained by Knecht \cite{Kn} for del Pezzo
surfaces of degree 2.  Also inspired by the methods of Hassett and
Tschinkel \cite{HaTs2}, Xu \cite{Xu} proved weak approximation for
del Pezzo surfaces of degree $\le 2$ that admit a model with a
certain nice property. For arbitrary del Pezzo surfaces of degree
$\le 3$ the question is still open.

As to the case where $k$ is not algebraically closed, particular
results for some special fields (e.g., the reals $\mathbb{R}$) have
been established with somewhat different methods (cf. \cite{CT},
\cite{Sch}, \cite{Du}).

In this paper, we study the case where the ground field $k$ is a
large field (cf. Definition \ref{defn2p3}) or a finite field. The
approach we follow was suggested by Colliot-Th\'{e}l\`{e}ne and
Hassett in 2005. It heavily relies on the method of Hassett and
Tschinkel \cite{HaTs1} as well as on some results of Koll\'{a}r and
Szab\'{o} (\cite{Kol2}, \cite{Kol3}, \cite{KoSz}).

Our main results are the following.

\begin{thm}\label{thm1p1}
With notation as in Problem$\;\ref{prob1p1}$, let $k$ be an infinite
algebraic extension of a finite field. Assume that
$X(K)\neq\emptyset$.

Then the $K$-variety $X$ satisfies weak approximation at places of
good reduction.
\end{thm}

\begin{thm}\label{thm1p2}
With notation as in Problem$\;\ref{prob1p1}$, let $k$ be a
nonarchimedean local field. Suppose that $\pi: \mathcal{X}\to C$ has
a section.

If for some point $p\in C(k)$, $R$-equivalence $($cf.
$\rm{Definition}\;\ref{defn2p5})$ on the rational points of the
fiber $\mathcal{X}_p$ is trivial, then there is a subset $W\subseteq
C(k)$ of rational points which is Zariski dense in $C$ such that the
$K$-variety $X$ satisfies weak approximation at places in $W$.
\end{thm}

Applying the methods used in the proofs, we also obtain the
following theorems.
\begin{thm}\label{thm1p3}
With notation as in Problem$\;\ref{prob1p1}$, let $k$ be a large
field. If $X(K)\neq\emptyset$, then $X(K)$ is Zariski dense in the
$K$-variety $X$.
\end{thm}

\begin{thm}\label{thm1p5}
Let $k=\mathbb{F}$ be a finite field, $C$ a smooth projective
$($geometrically connected$)$ curve over $k$, $K=k(C)$ its function
field, and $\mathcal{X}\subseteq{\mathbb{P}^{N}_k\times_kC}$ a
smooth closed subvariety such that the canonical projection $\pi:
\mathcal{X}\to C$ is a surjective morphism whose generic fiber is a
smooth cubic hypersurface $X\subseteq\mathbb{P}^N_K$ of dimension
$\ge 2$. Suppose that $X(K)\neq\emptyset$.

Then for any finite set $P$ of closed points of good reduction on
$C$ with
\[
|\kappa(p)|>10\,,\quad\forall\; p\in P\,,
\] the specialization map
\begin{equation}\label{eq1p1}
X(K)\longrightarrow \prod_{p\in P}\mathcal{X}_{p}(\kappa(p))
\end{equation}is surjective. In particular, if $|k|>10$, then question
$(\rm{SP})$ has a positive answer in this situation.
\end{thm}

In proving the above theorem, we use some ideas of Swinnerton-Dyer,
who proved an analogous result for cubic surfaces over number fields
(cf. \cite{SD}, Theorem 5). In particular, his method of descent
from a tower of quadratic extensions of the ground field finds a
significant application in our proof. Yet another important
ingredient is a result of Koll\'ar on the geometry of cubic
hypersurfaces over finite fields (\cite{Kol4}, Lemma$\;$9.4).

As a consequence, we obtain a weak approximation theorem for cubic
surfaces, which is analogous to a result of Swinnerton-Dyer
(\cite{SD}, Corollary to Theorem 5).

\begin{thm}\label{thm1p6}
Let $k=\mathbb{F}_q$ be a finite field of characteristic not
dividing $6$ and of cardinality $q>47$, $C$ a smooth projective
$($geometrically connected$)$ curve over $k$, $K=k(C)$ its function
field, and $\mathcal{X}\subseteq{\mathbb{P}^{3}_k\times_kC}$ a
smooth closed subvariety such that the canonical projection $\pi:
\mathcal{X}\to C$ is a surjective morphism whose generic fiber is a
smooth cubic surface $X\subseteq\mathbb{P}^3_K$ with
$X(K)\neq\emptyset$. Let $v$ be a place of good reduction
corresponding to a closed point $p\in C$.

Then $X(K)$ is $v$-adically dense in $X(K_v)$. More precisely, for
any $x\in \mathcal{X}_p(\kappa(p))$ there is an $R$-equivalence
class $\mathscr{C}$ on $X(K)$ such that $\mathscr{C}$ is
$v$-adically dense in the inverse image of $x$ under the
specialization map $X(K_v)\to \mathcal{X}_{p}(\kappa(p))$.
\end{thm}

\section{Review of Basic Notions}

\subsection{Separably rational connectedness}

\begin{defn}[\cite{Deb}, Definition$\;$4.5]
Let $X$ be a smooth variety and $r\ge 0$. An \empha{$r$-free
rational curve}\index{free curve} on $X$ is a nonconstant morphism
$f: \mathbb{P}^1\to X$ such that
\[
f^*T_X\cong \mathscr{O}_{\mathbb{P}^1}(a_1)\oplus\cdots\oplus
\mathscr{O}_{\mathbb{P}^1}(a_n)
\]with $a_1\ge \cdots\ge a_n\ge r$. (By a theorem of Grothendieck (\cite{Har}, p.384, Exercise$\;$2.6), every vector
bundle over $\mathbb{P}^1$ is isomorphic to a direct sum
$\oplus\mathscr{O}_{\mathbb{P}^1}(a_i)$ and the sequence
$a_1\ge\cdots\ge a_n$ is uniquely determined by the vector bundle.)
We shall say \empha{free} (resp. \empha{very free}) instead of
0-free (resp. 1-free).
\end{defn}

\begin{defn}[\cite{Kol}, IV.3.2]\label{defn2p1}
Let $k$ be a field and $\Bar{k}$ an algebraic closure of $k$. A
smooth proper $k$-variety $X$ is called \empha{separably rationally
connected} if it has the following property:

(1) There is a separated integral algebraic $k$-scheme $Z$ and a
morphism $g: \mathbb{P}^1\times_kZ\to X$ such that the morphism
\[
g^{(2)}:\;\;
\mathbb{P}^1\times_k\mathbb{P}^1\times_kZ\longrightarrow X\times_k
X\,;\;\;(t_1\,,\,t_2\,,\,z)\mapsto (g(t_1\,,\,z)\,,\,g(t_2\,,\,z))
\]is dominant and smooth at the generic point.

When $\mathrm{dim}\,X>0$, this is equivalent to the following two
properties:

(2) There is a very free rational curve $\mathbb{P}^1_{\Bar{k}}\to
X_{\Bar{k}}$ defined over $\Bar{k}$.

(3) For any finite collection of closed points $x_1,\dotsc, x_m\in
X_{\Bar{k}}$ and any tangent directions $\xi_i$ at the points $x_i$,
there is a very free curve $f: \mathbb{P}^1_{\Bar{k}}\to
X_{\Bar{k}}$ defined over $\Bar{k}$ such that $f$ is an unramified
morphism with image containing all the $x_i$ as smooth points and
having tangent directions $\xi_i$ at these points.
\end{defn}

In property (3) of the above definition, one can take $f$ to be a
closed embedding if $\dim\, X\ge 3$. For a proof of the above
equivalences, see \cite{Kol} IV.3.7 and 3.9, and \cite{Deb2},
Th\'eor\`eme$\;$2.2.

\begin{examples}\

(1) A smooth $k$-variety $X$ is separably rationally connected if it
is geometrically separably unirational (i.e., over $\Bar{k}$ there
is a separable dominant rational map from some
$\mathbb{P}^n_{\Bar{k}}$ to $ X_{\Bar k}$).

(2) In characteristic 0, smooth Fano varieties are separably
rationally connected (Campana, \cite{Cam}). In particular, a smooth
hypersurface in $\mathbb{P}^n$ of degree $d\le n$ is separably
rationally connected.

(3) In characteristic $p>0$, smooth cubic hypersurfaces of dimension
$\ge 2$ are separably rationally connected. (Smooth cubic surfaces
are geometrically rational. For the case of dimension $\ge 3$, see
\cite{Con} Remarque$\;$2.34.)
\end{examples}

\subsection{Large fields and $R$-equivalence}

\begin{defn}[\cite{FJ}, Chapter 10 and \cite{Pop}]\label{defn2p3}
A field $k$ is said to be

(1) \empha{pseudo-algebraically closed} (PAC), if every
(geometrically integral) $k$-variety  has a rational point;

(2) \empha{large}\index{large field}, if for every (geometrically
integral) $k$-variety $X$ having a smooth $k$-point, $X(k)$ is
Zariski dense in $X$.
\end{defn}

\begin{examples}
\

(1) PAC fields are large fields. Algebraic extensions of PAC (resp.
large) fields are PAC (resp. large). (cf. \cite{FJ}, Chapter 10 and
\cite{Pop}, Proposition$\;$1.2.)

(2) Infinite algebraic extensions of finite fields are PAC. (cf.
\cite{FJ}, Chapter 10.)

(3) The class of large fields contains all local fields (i.e.,
finite extensions of $\mathbb{R}$, $\mathbb{Q}_p$ or the field of
Laurent series $\mathbb{F}_p((t))$ over the finite field
$\mathbb{F}_p$). More generally, a field that is complete with
respect to a nontrivial absolute value is large.

(4) All real closed fields and all $p$-adic closed fields (i.e.,
fields with absolute Galois group isomorphic to that of $\mathbb{R}$
or that of a $p$-adic field) are large (cf. \cite{PrRo},
Theorem$\;$7.8).
\end{examples}

\begin{defn}[Manin]\label{defn2p5}
Let $X$ be an algebraic scheme over a field $k$. Two rational points
$x,\,y\in X(k)$ are said to be \empha{directly $R$-equivalent} if
there is a $k$-rational map $f: \mathbb{P}^1_k\dashrightarrow X$
such that $f(0)=x$ and $f(\infty)=y$. The \empha{$R$-equivalence
relation} on $X(k)$ is the equivalence relation generated by direct
$R$-equivalence. We will write $X(k)/R$ for the set of
$R$-equivalence classes on $X(k)$. When $X(k)/R$ has only one
element, we write $X(k)/R=1$ and say that $R$-equivalence is trivial
on $X$.
\end{defn}

$R$-equivalence is closely related to separably rational
connectedness. In \cite{Kol2}, Koll\'{a}r proved the finiteness of
$R$-equivalence classes on separably rationally connected varieties
over local fields. Later he generalized the main theorem in
\cite{Kol2} as follows.

\begin{thm}[{Koll\'{a}r, \cite{Kol3}}, Theorem$\;$23]\label{thm2p6}
Let $k$ be a large field and let $X$ be a smooth, projective,
separably rationally connected $k$-variety of dimension $\ge 3$. Let
$n\ge 1$ and $x_1,\dotsc, x_n\in X(k)$. Then the following two
conditions are equivalent:

$(\mathrm{i})$ There is a closed embedding $f: \mathbb{P}^1\to X$
that is $(n+1)$-free such that $x_1,\dotsc, x_n\in
f(\mathbb{P}^1(k))$.

$(\mathrm{ii})$ All the points $x_i\in X(k)$ are in the same
$R$-equivalence class.
\end{thm}

We remark that in general it may happen that there is no smooth
rational curve  passing through a given finite set of rational
points on a variety $X$ as in the above theorem. A sufficient
condition for this is that the ground field $k$ is PAC (cf.
\cite{KoSz}, Theorem$\;$19 and \cite{Wi}, Corollaire$\;$3.2). If
$k=\mathbb{F}$ is a finite field, we have the following result.

\begin{thm}[{Koll\'{a}r and Szab\'{o}, \cite{KoSz}}, Theorems$\;$2 and
7]\label{thm2p7} \

$(\mathrm{i})$ There is a function $\Phi : \mathbb{N}^3\to
\mathbb{N}$ having the following property: Let
$X\subseteq\mathbb{P}^N$ be a smooth, projective, separably
rationally connected variety over a finite field $\mathbb{F}$ of
dimension $\ge 3$, and let $S\subseteq X$ be a zero-dimensional
smooth subscheme. If $|\mathbb{F}|>\Phi(\deg\, X,\,\dim X,\,\deg\,
S)$, then there is a closed embedding $f: \mathbb{P}^1\to X$ that is
$(\deg\, S+1)$-free with image containing $S$.

$(\mathrm{ii})$ There is a function $\Psi:
\mathbb{N}^2\to\mathbb{N}$ having the following property: Let
$X\subseteq\mathbb{P}^N$ be a smooth, projective, separably
rationally connected variety over a finite field $\mathbb{F}$. If
$|\mathbb{F}|>\Psi(\deg\, X,\,\dim X)$, then $X(\mathbb{F})/R=1$.
\end{thm}

In the above two theorems the hypothesis that $\dim X\ge 3$ is not
essential. It permits us to take $f$ to be a closed embedding. If
$\dim X=2$, we can apply the theorem to $X\times_k\mathbb{P}^1$ to
obtain a rational curve passing through the given points.

\begin{remark}\label{remarkAdded1}
In Theorem$\;$\ref{thm2p7} (ii), if one restricts to the case where
$X\subseteq\mathbb{P}^N$ is a smooth cubic hypersurface of dimension
$\ge 2$, then \cite{Kol4}, Theorem$\;$1.1 gives an even better lower
bound for $|\mathbb{F}|$: one has $X(\mathbb{F})/R=1$ whenever
$|\mathbb{F}|\ge 11$.\footnote{The lower bound $|\mathbb{F}|\ge 8$
originally given in \cite{Kol4} should be $|\mathbb{F}|\ge 11$. See
also the footnote attached to Lemma$\;$\ref{lemmaAdded1}.}
\end{remark}

\section{Weak approximation over function fields}
To establish weak approximation we use the method initiated by
Hassett and Tschinkel in \cite{HaTs1}. In this section we briefly
recall their method.

\subsection{Geometric interpretation}
Let $K=k(C)$ be the function field of a smooth projective curve and
let $X$ be any smooth proper variety over $K$. A result of Nagata
(\cite{Nag}) guarantees that $X$ always admits a \empha{proper
model}, i.e., an algebraic space $\pi: \mathcal{X}\to C$ that is
flat and proper over $C$ with generic fiber $X$. Rational points of
$X$ correspond bijectively in a natural way to sections of $\pi:
\mathcal{X}\to C$.

In this paper we will assume for simplicity that \emph{$X$ is
projective and admits a smooth projective model} $\pi:
\mathcal{X}\to C$ as in Problem$\;$\ref{prob1p1}.

 For each closed point $p\in C$, let
$\mathcal{X}_p=\mathcal{X}\times_C\mathrm{Spec}(\kappa(p))$ be the
fiber over $p$. Let $\widehat{\mathscr{O}}_{C,\,p}$ be the completed
local ring at $p$, with maximal ideal
$\widehat{\mathfrak{m}}_{C,\,p}$. If $v$ denotes the place of $K$
corresponding to $p$, then $K_v$ is the field of fractions of
$\widehat{\mathscr{O}}_{C,\,p}$. Let
$\widehat{C}_p=\mathrm{Spec}(\widehat{\mathscr{O}}_{C,\,p})$,
$\mathcal{X}_{\widehat{C}_p}=\mathcal{X}\times_C\widehat{C}_p$ and
$\widehat{\pi}_p: \mathcal{X}_{\widehat{C}_p}\to \widehat{C}_p$ the
natural projection. Sections of $\widehat{\pi}_p$ correspond
bijectively to $K_v$-valued points of $X$. Suppose a point $M_v\in
X(K_v)$ corresponds to a section $\widehat{s}_p:
\widehat{C}_p\to\mathcal{X}_{\widehat{C}_p}$. Then basic $v$-adic
open neighborhoods of $M_v$ consist of those sections of
$\widehat{\pi}_p$ which agree with
\[
\widehat{C}_{p,\,N}:=\mathrm{Spec}(\widehat{\mathscr{O}}_{C,\,p}/\widehat{\mathfrak{m}}_{C,\,p}^{N+1})\hookrightarrow\widehat{C}_p
\overset{\widehat{s}_p}{\longrightarrow}\mathcal{X}_{\widehat{C}_p}
\]when restricted to $\widehat{C}_{p,\,N}\subseteq\widehat{C}_p$, for some
$N\in\mathbb{N}$.

For $N\in\mathbb{N}$, we say a morphism
\[
j: \;\;
\widehat{C}_{p,\,N}=\mathrm{Spec}(\widehat{\mathscr{O}}_{C,\,p}/\widehat{\mathfrak{m}}_{C,\,p}^{N+1})
\longrightarrow\mathcal{X}_{\widehat{C}_p}=\mathcal{X}\times_C\mathrm{Spec}(\widehat{\mathscr{O}}_{C,\,p})
\]is an \empha{$N$-jet}\index{jet} at the closed point $p\in C$ if $\widehat{\pi}_p\circ j$ coincides with the
natural inclusion
$\widehat{C}_{p,\,N}\hookrightarrow\widehat{C}_p=\mathrm{Spec}(\widehat{\mathscr{O}}_{C,\,p})$.
A section of $\widehat{\pi}_p$ determines naturally an $N$-jet.

Thus weak approximation over function fields has a geometric
reformulation as follows:  for any finite set of closed points
$\{p_i\}$ in $C$, any sections $\widehat{s}_i$ of
$\widehat{\pi}_{p_i}$ and any $N\in\mathbb{N}$, there exists a
section $s$ of $\pi$ such that for each $i$, the $N$-jets determined
by $s$ and $\widehat{s}_i$ coincide.

\subsection{Iterated blowup and weak approximation}
We keep notation as above. To develop the method of Hassett and
Tschinkel over an arbitrary ground field $k$, we assume that places
in the consideration below are given by \emph{separable points}
(i.e. points with residue fields separable over $k$).

Let $\{\,p_i\,\}\subseteq C$ be a finite number of separable closed
points of good reduction and let $J=\{\,j_i\,\}$ be $N$-jets at
these points given by formal sections $\widehat{s}_i$ of
$\widehat{\pi}_{p_i}: \mathcal{X}_{\widehat{C}_{p_i}}\to
\widehat{C}_{p_i}$. By considering the \empha{iterated blowup}
associated to $J$, whose construction we will recall below, Hassett
and Tschinkel reduce finding a section passing through given points
with prescribed jet data to simply finding a section through given
points on the iterated blowup.

Here is the construction of the iterated blowup $\beta(J):
\mathcal{X}(J)\to \mathcal{X}$ associated to $J$. It is obtained by
 the following sequence of blowups:
\[
\mathcal{X}(J)=\mathcal{X}_{(N)}\to \mathcal{X}_{(N-1)}\to\cdots\to
\mathcal{X}_{(1)}\to\mathcal{X}
\]where $\mathcal{X}_{(1)}\to \mathcal{X}$ is the blowup of $\mathcal{X}$ at
$\{\,\widehat{s}_i(p_i)\,\}$ and for each $n=1,\,\dotsc,\, N-1$,
$\mathcal{X}_{(n+1)}\to\mathcal{X}_{(n)}$ is the blowup at the
points where the proper transform of $\widehat{s}_i$ meets the fiber
over $p_i$.

Write $d=\dim X$. For each $i$ the fiber $\mathcal{X}(J)_{p_i}$
decomposes into irreducible components:
\[
\mathcal{X}(J)_{p_i}=E_{i,\,0}\cup \cdots\cup E_{i,\,N}
\]where

(i) $E_{i,\,0}$ is the proper transform of $\mathcal{X}_{p_i}$,
isomorphic to the blowup of $\mathcal{X}_{p_i}$ at the
$\kappa(p_i)$-point $r_{i,\,0}:=\widehat{s}_i(p_i)$;

(ii) $E_{i,\,n}$ for $n=1,\dotsc, N-1$ is the blowup of
$\mathbb{P}^{d}_{\kappa(p_i)}$ at the $\kappa(p_i)$-point
$r_{i,\,n}$ where the proper transform of $\widehat{s}_i$ meets the
fiber over $p_i$ of the $n$-th blowup;

(iii) $E_{i,\,N}\cong\mathbb{P}^{d}_{\kappa(p_i)}$.

Let $r_{i}\in  E_{i,\,N}\setminus E_{i,\,N-1}$ denote the point
where the proper transform of $\widehat{s}_i$ meets $E_{i,\,N}$.
Sections $s'$ of $\mathcal{X}(J)\to C$ with $s'(p_i)=r_i$ yield
sections of $\pi: \mathcal{X}\to C$ with $N$-jets $j_i$ at the
points $p_i$. So we have the following criterion.

\begin{prop}[\cite{HaTs1}, Proposition$\;$11]\label{prop3p1}
Let $X,\,C$ and $\pi: \mathcal{X}\to C$ be as above. Consider all
data $(\{\,p_i\,\},\,J=\{\,j_i\,\},\, \{\,r_i\,\})$ consisting of a
finite collection of separable closed points $\{\,p_i\,\}\subseteq
C$ of good reduction, $N$-jets $\{\,j_i\,\}$ at these points and
$\kappa(p_i)$-points $\{\,r_i\,\}$ on the corresponding iterated
blowup $\mathcal{X}(J)$ with $r_i\in E_{i,\,N}\setminus
E_{i,\,N-1}$.

If for every datum there exists a section $s'$ of $\mathcal{X}(J)\to
C$ with $s'(p_i)=r_i$ for all $i$, then $X$ satisfies weak
approximation at places given by separable closed points of good
reduction.
\end{prop}

\subsection{Deformations of combs}

The main steps in constructing sections through prescribed points
involve techniques for smoothing combs.

\begin{defn}[\cite{HaTs1}, Definition$\;$18]
Let $k$ be any field and let $\Bar{k}$ be an algebraic closure of
$k$. A \empha{comb with $m$ $($reducible$)$ teeth} over $k$ is a
projective curve $C^*$ over $k$ with a distinguished irreducible
component $D\subseteq C^*$ defined over $k$ such that the following
conditions hold:

(1) $D$ is a smooth (geometrically connected) curve over $k$, called
the \empha{handle} of the comb;

(2) $\overline{C}^*=C^*\times_k\Bar{k}$ is the union of
$\overline{D}=D\times_k\Bar{k}$ and $m$ other curves
$\overline{T}_1,\dotsc, \overline{T}_m$ which are called the
\empha{teeth} of the comb;

(3) each $\overline{T}_j$ is a chain of $\mathbb{P}^1$'s;

(4) the $\overline{T}_j$ are disjoint with one another and each of
them meets $\overline{D}$ transversally at a single point.
\end{defn}

Let $X,\,C$, $\pi: \mathcal{X}\to C$, $\{\,p_i\,\}$ and so on be as
in Proposition$\;$\ref{prop3p1}. Recall that we are assuming that
the total space $\mathcal{X}$ is smooth projective and that the
$p_i$ are separable points. Since $X$ satisfies weak approximation
if and only if $X\times_K\mathbb{P}^1_K$ does, without loss of
generality we may assume that $d=\dim X\ge 3$.

Suppose we are given a section $\sigma: C\to\mathcal{X}(J)$. Let
$I'$ (resp. $I''$) denote the set of indices with
$q_i:=\sigma(p_i)\neq r_i$ (resp. $q_i=r_i$).

\

Working over an algebraically closed ground field $k$, the most
technical part of \cite{HaTs1} shows that a section
$C\to\mathcal{X}(J)$ passing through $r_i$ can be found so long as
the following two tasks (A) and (B) may be accomplished:

(A). In the case $N=0$, find a comb $C^*$ with handle $D=\sigma(C)$
and smooth teeth $T_1,\dotsc, T_m$ and a closed embedding $f: C^*\to
\mathcal{X}$ such that the following properties hold:

$(\mathrm{i})$ the teeth $T_j$ are very free rational curves
contained in distinct general fibers of $\pi$;

$(\mathrm{ii})$ for each $i\in I'$, there is a tooth $T_i$ with
image containing $r_i$ as a smooth point;

$(\mathrm{iii})$ if $\mathfrak{r}$ denotes the set of points on
$C^*$ mapping to the $r_i,\,i\in I=I'\cup I''$ and
$N_f=N_{C^*/\mathcal{X}}$ denotes the normal bundle of $C^*$ in
$\mathcal{X}$, then the restriction of
$N_{f}\otimes\mathscr{O}_{C^*}(-\mathfrak{r})$ to every irreducible
component of $C^*$ is generated by global sections and has no higher
cohomology.

(B). In the case $N\ge 1$, assume $q_i\in E_{i,\,N}\setminus
E_{i,\,N-1}$. Find a comb $C^*$ with handle $D=\sigma(C)$ and
reducible teeth $T_1,\dotsc, T_m$ and a closed embedding $f: C^*\to
\mathcal{X}(J)$ having the following properties:

$(\mathrm{i})$ for each $i\in I'$, there is a tooth $T_i$ mapped
into $\mathcal{X}(J)_{p_i}$ and containing $r_i$;

$(\mathrm{ii})$ for each $i\in I'$, $f(C^*)$ contains $r_i$ as a
smooth point, so there is a unique irreducible component
$T_{i,\,N}\subseteq C^*$ with image containing $r_i$;

$(\mathrm{iii})$ the remaining teeth $T_j,\, j=|I'|+1,\dotsc, m$ are
very free rational curves contained in distinct general fibers of
$\mathcal{X}(J)\to C$;

$(\mathrm{iv})$ if $\mathfrak{r}$ denotes the set of points on $C^*$
mapping to the $r_i$, $i\in I=I'\cup I''$ and
$N_f=N_{C^*/\mathcal{X}(J)}$ denotes the normal bundle of $C^*$ in
$\mathcal{X}(J)$, then the restriction of
$N_f\otimes\mathscr{O}_{C^*}(-\mathfrak{r})$ to every irreducible
component of $C^*$ is generated by global sections and has no higher
cohomology.

In \cite{HaTs1}, (A) and (B) are done in Lemmas$\;$25 and 26,
assuming of course the generic fiber $X$ is separably rationally
connected. The desired section is a smooth deformed curve of the
comb $C^*$, which is obtained by applying \cite{HaTs1},
Proposition$\;$24.

The key condition here is the vanishing of the higher cohomology of
the twisted normal sheaf
$N_f\otimes\mathscr{O}_{C^*}(-\mathfrak{r})$. We consider this
condition on each irreducible component of $C^*$. For a component in
a tooth $T_j$, the condition is easier to verify since such a curve
is usually a very free curve. The difficult part is the restriction
to the handle $D=\sigma(C)$. Fortunately, the sheaf
$N_f\otimes\mathscr{O}_D$ has a nice interpretation as a kind of
elementary transform of $N_{\sigma}$, where $N_\sigma$ denotes the
normal sheaf of $D$ in $\mathcal{X}(J)$. Let $D^c$ be the closure of
$C^*\setminus D$ in $C^*$ and $S=D\cap D^c$ the locus where the
handle $D$ meets the teeth. For each $q_i\in S$, let $\xi_i\subseteq
N_\sigma\otimes\kappa(q_i)$ denote the subspace determined by the
tangent direction $T_{D^c,\,q_i}$. It turns out that
$N_f\otimes\mathscr{O}_D$ is the so-called \empha{sheaf of rational
sections of $N_{\sigma}$ having at most a simple pole at each $q_i$
in the direction $\xi_i$ and regular elsewhere} (cf. \cite{GHS},
Lemma$\;$2.6).

A key lemma needed is the following.

\begin{lemma}[\cite{GHS}, Lemma$\;$2.5]\label{lemma3p3}
Let $C$ be a smooth curve over an algebraically closed field,
$\mathscr{V}$ a vector bundle on $C$ and let $l\ge 1$ be a positive
integer. Then for  sufficiently large $m$ there exist general points
$q_1,\dotsc, q_m\in C$ and one-dimensional subspaces $\xi_i\subseteq
\mathscr{V}\otimes\kappa(q_i)$ such that for any $l$ points
$w_1,\dotsc, w_l\in C$ the sheaf $\mathscr{V}\,'(-w_1-\cdots-w_l)$
is generated by global sections and
\[
H^1(C,\,\mathscr{V}\,'(-w_1-\cdots-w_l))=0\,,
\]where $\mathscr{V}\,'$ denotes the sheaf of rational sections of
$\mathscr{V}$ that have at most a simple pole at each $q_i$ in the
direction $\xi_i$ and regular elsewhere.
\end{lemma}

Note that the smallest value for $m$ depends on the number $l$ and
the vector bundle $\mathscr{V}$.

Lemma$\;$\ref{lemma3p3} applied to the normal bundle
$\mathscr{V}=N_\sigma$ on $D$ shows that the vanishing of higher
cohomology required in the tasks (A) and (B) may be guaranteed by
adding sufficiently many teeth to the comb $C^*$.

\

Now we want to deal with the case where $k$ is not necessarily
algebraically closed and we want to produce a deformed curve
\emph{defined over $k$}. We face two problems.

(D.1) We should be able to construct the comb $C^*$ over $k$.

First consider task (A). Suppose that we want to add teeth in the
fibers over separable closed points $p_i\in C$.  If we are able to
find in each fiber a very free curve $T_i$ defined over
$\kappa(p_i)$, then the union $\cup_iT_i\cup\sigma(C)$ considered as
a $k$-scheme gives us a comb defined over $k$. Unlike the case of
algebraically closed ground field, the existence of a very free
curve defined over the ground field in a separably rationally
connected fiber is not a priori evident, and if we moreover require
the very free curve to pass through two prescribed points, this may
be even impossible. In the next section we shall use the results of
Koll\'{a}r and Szab\'{o} on the triviality of $R$-equivalence to
solve this problem.

The situation for (B) is essentially the same. The difference is
that the fiber $\mathcal{X}(J)_{p_i}$ becomes a bit more
complicated. However, by a careful inspection on the construction of
reducible teeth in \cite{HaTs1}, we see that for each reducible
tooth $T_i$ all but one of its components $T_{i,\,0},\dotsc,
T_{i,\,N}$ are just lines connecting two $\kappa(p_i)$-points and
the last one $T_{i,\,0}$ is a very free curve passing through one
$\kappa(p_i)$-point in the piece $E_{i,\,0}$, which is the blowup of
$\mathcal{X}_{p_i}$ at a $\kappa(p_i)$-point. Separably rational
connectedness and the set of $R$-equivalence classes are both
birationally invariant for smooth projective varieties (cf.
\cite{Kol}, IV.3.3 and \cite{CTSan}, p.195, Proposition$\;$10). So
if a property determined by separably rational connectedness or
$R$-equivalence holds for $\mathcal{X}_{p_i}$, then the same should
hold for $E_{i,\,0}$. In particular, if the same thing can be done
in $\mathcal{X}_{p_i}$, then the curve $T_{i,\,0}$ in $E_{i,\,0}$
may be constructed over $\kappa(p_i)$. So it makes no essential
difference in fulfilling tasks (A) and (B).

(D.2) We should be able to find deformations of $C^*$ that are
defined over $k$.

Having constructed a comb $C^*$ defined over $k$ with good
properties as stated in (A) or (B), by applying the same argument as
in \cite{HaTs1} we see that the parameter space $M$ of curves in
$\mathcal{X}(J)$ passing through the separable points $\{\,r_i\,\}$
is smooth at the $k$-point $[C^*]$ and there are nearby deformed
curves defined over $\Bar{k}$ that are smooth. The final step is to
show the existence of a $k$-rational point in the parameter space.
When the field $k$ is large or finite, we have some tools to do this
as we shall now see.

\section{Proofs of Results over Large Fields}

\begin{prop}\label{prop4p1}
With notation as in Problem$\;\ref{prob1p1}$, let $k$ be a large
field. Consider question $(\rm{SP})$. Suppose that the closed points
$p_1,\dotsc, p_n$ are separable and that the point $r_i\in
\mathcal{X}_i(\kappa(p_i))$ lies in the same $R$-equivalence class
as $q_i=\sigma(p_i)\in \mathcal{X}_i(\kappa(p_i))$ for every
$i=1,\dotsc, n$.

Then there exists a section $s: C\to \mathcal{X}$ of $\pi$ such that
$s(p_i)=r_i$ for $i=1,\dotsc, n$.
\end{prop}
\begin{proof}
In each fiber $\mathcal{X}_i=\mathcal{X}_{p_i}$, the two rational
points $r_i$ and $q_i=\sigma(p_i)$ are $R$-equivalent by hypothesis.
Thus Koll\'{a}r's theorem (Theorem$\;$\ref{thm2p6}) implies that
there is a very free curve $f_i: \mathbb{P}^1\to \mathcal{X}_i$
defined over $\kappa(p_i)$ such that $q_i\,,\,r_i\in
f_i(\mathbb{P}^1(\kappa(p_i)))$. On the other hand, the smooth curve
$C$ has sufficiently many separable closed points. So we can pick as
many good points as we need and construct teeth in the corresponding
fibers to get a comb $C^*$ defined over $k$ satisfying all the
properties in task (A). The comb $C^*$ then has smooth curves as
nearby deformations, and since $[C^*]$ is a smooth $k$-point in the
parameter space, the property of large fields implies that there
exists a deformed smooth curve that is defined over $k$. The
problems (D.1) and (D.2) described in $\S$3 may be thus solved. We
then get a section $s$ passing through the points $r_i$.
\end{proof}

\begin{proof}[Proof of Theorem$\;\ref{thm1p1}$]
Let $U\subseteq C$ be the open subset consisting of points of good
reduction. (Such an open subset exists by \cite{Kol}, IV.3.11.) Note
that $k$ is perfect so that every point is separable. Let
$p_1,\dotsc, p_n$ be closed points lying in $U$.

The theorem of Koll\'{a}r and Szab\'{o} (Theorem$\;$\ref{thm2p7}
(ii)) implies that $R$-equivalence on a smooth separably rationally
connected variety over an infinite algebraic extension of a finite
field is trivial. So we have $\mathcal{X}_{p_i}(\kappa(p_i))/R=1$. A
similar argument as in the proof of Proposition$\;$\ref{prop4p1}
shows that a comb $C^*$ defined over $k$ may be constructed
satisfying all the required properties in task (A) or (B). Since $k$
is a large field, the problem (D.2) may also be solved. This means
that the comb $C^*$ admits smooth deformed curves defined over $k$,
which give rise to sections on the iterated blowup with the required
property. Weak approximation at places in $U$ is thus proved for the
generic fiber $X$.
\end{proof}

\begin{proof}[Proof of Theorem$\;\ref{thm1p2}$]
By a result of Koll\'{a}r (\cite{Kol3}, Theorem$\;$3), the function
\[
C(k)\longrightarrow \mathbb{N}\,,\quad\;\; c\mapsto
|\mathcal{X}_c(k)/R|
\]is upper semi-continuous for the $k$-topology. So there is a nonempty open subset $W$ of $C(k)$ (for the
$k$-topology) such that $\mathcal{X}_{c}(k)/R=1$ for all $c\in W$.
For any nonempty Zariski open subset $U$ of $C$, $U(k)$ is dense in
$C(k)$ for the $k$-topology. Hence $W\cap U(k)\neq \emptyset$,
showing that $W$ is Zariski dense in $C$. The same argument as in
the proof of Theorem$\;$\ref{thm1p1} proves weak approximation at
places  in $W$.
\end{proof}

\begin{remark}
In Theorem$\;$\ref{thm1p2}, the hypothesis on the existence of a
point $p\in C(k)$ with $\mathcal{X}_p(k)/R=1$ is satisfied if there
is a point $p\in C(k)$ such that $\mathcal{X}_p$ has a smooth
separably rationally connected reduction over the residue field
$\mathbb{F}$ of $k$ and if the residue field $\mathbb{F}$ has
sufficiently large cardinality (cf. \cite{KoSz}, Theorem$\;$8).
\end{remark}

The following lemma is an easy consequence of a theorem of
Koll\'{a}r.

\begin{lemma}\label{lemma4p3}
Let $X$ be a smooth, projective, separably rationally connected
variety over a large field $k$. Suppose that $X(k)\neq \emptyset$.
Then for every $x\in X(k)$, there is a subset $W\subseteq X(k)$ that
is Zariski dense in $X$ with the following property: for every $y\in
W$ there is a very free curve $f_{xy}: \mathbb{P}^1\to X$ defined
over $k$ such that $x,\,y\in f_{xy}(\mathbb{P}^1(k))$.
\end{lemma}
\begin{proof}
In \cite{Kol2}, Theorem$\;$1.4, Koll\'{a}r proved that there exists
a very free curve $f_x: \mathbb{P}^1\to X$ defined over $k$ sending
a rational point  $0\in \mathbb{P}^1$ to $x$. Let
$\mathrm{Hom}(\mathbb{P}^1,\,X\,; 0\mapsto x)$ be the Hom-scheme
parametrizing morphisms $\mathbb{P}^1\to X$ sending $0$ to $x$. By
\cite{Deb}, p.91, Proposition$\;$4.8, the evaluation map
\[
e: \quad \mathbb{P}^1\times\mathrm{Hom}(\mathbb{P}^1,\,X\,; 0\mapsto
x)\longrightarrow X\,,\quad (t,\, g)\mapsto g(t)
\]is smooth at any rational point $(t,\,[f_x])$ with $t\neq 0$. In
particular, $e$ is dominant. Therefore, for every nonempty open
subset $U$ of $X$, we can find a very free curve $f: \mathbb{P}^1\to
X$ defined over $k$ sending a rational point of $\mathbb{P}^1$ into
$U$ with $f(0)=x$.
\end{proof}

\begin{proof}[Proof of Theorem$\;\ref{thm1p3}$]
Let $\sigma : C\to\mathcal{X}$ be a section of $\pi: \mathcal{X}\to
C$. Let $Z$ be a proper closed subset of $X$. Its closure in
$\mathcal{X}$ is a proper closed subset $\mathcal{Z}$ of
$\mathcal{X}$. There is a nonempty open subset $U$ of $C$ such that
for every $c\in U$, the fiber $\mathcal{X}_c$ is smooth separably
rationally connected and $\mathcal{X}_c\cap \mathcal{Z}\neq
\mathcal{X}_c$. Choose a separable closed point $p\in U$. By
Lemma$\;$\ref{lemma4p3}, there is a point $r\in
\mathcal{X}_p(\kappa(p))\setminus \mathcal{Z}$ which is directly
linked to $q=\sigma(p)$ by a very free rational curve defined over
$\kappa(p)$. Proposition$\;$\ref{prop4p1} implies that there is a
section $s: C\to \mathcal{X}$ such that $s(p)=r$. This section gives
a $K$-rational point of $X$ which does not lie in $Z$. Since $Z$ is
arbitrary, it follows that $X(K)$ is Zariski dense in $X$.
\end{proof}

\section{Proofs of Results over Finite Fields}
\subsection{Surjectivity of specialization}
\begin{prop}\label{prop1p4}
With notation as in Problem$\;\ref{prob1p1}$, let $k=\mathbb{F}$ be
a finite field. Let $S\subseteq C$ be the closed subset consisting
of points of bad reduction for the model $\pi: \mathcal{X}\to C$.
Assume that $\pi$ has a section $\sigma: C\to \mathcal{X}$. We fix
closed embeddings
\[
C\hookrightarrow\mathbb{P}^3_k\,,\;\;\;
\mathcal{X}\hookrightarrow\mathbb{P}^N_k\times_kC\hookrightarrow\mathbb{P}^N_k\times_k\mathbb{P}^3_k\hookrightarrow
\mathbb{P}^M_k\] where
$\mathbb{P}^N_k\times_k\mathbb{P}^3_k\hookrightarrow \mathbb{P}^M_k$
is the Segre embedding. We denote by $\deg\,C$ the degree of $C$ in
$\mathbb{P}^3_k$, $\deg\,\sigma(C)$ and $\deg\,\mathcal{X}$ the
degrees of $\sigma(C)$ and $\mathcal{X}$ in $\mathbb{P}^M_k$, and
$\deg\,X$ the degree of $X$ with respect to the embedding
\[
X=\mathcal{X}\times_CK\;\hookrightarrow\;
\mathbb{P}^N_K=(\mathbb{P}^N_k\times_kC)\times_CK\,.
\]

There is a lower bound
\[
B=B(N_{\sigma},\,|S|\,, \deg_CP ,\,\deg\,
C,\,\deg\,\sigma(C)\,,\,\deg\, X,\,\dim X)
\] for
the cardinality $|\mathbb{F}|$, depending on the normal bundle
$N_{\sigma}$ of the section $\sigma$, the cardinality $|S|$, the
number $\deg_CP=\sum_{p\in P}[\kappa(p):k]$ and the geometric
invariants $\deg\,C,\,\deg\,\sigma(C),\,\deg\, X$ and $\dim\, X$,
such that the answer to question $(\rm{SP})$ is yes whenever
$|\mathbb{F}|>B$.
\end{prop}
\begin{proof}
According to the method we described in $\S$3, we want to first
construct a comb $C^*$ such that the conditions in task (A) hold.
When $|\mathbb{F}|$ is greater than a lower bound in terms of
$\deg\,X$ and $\dim X$, by the theorem of Koll\'{a}r and Szab\'{o}
(Theorem$\;$\ref{thm2p7} (i)), in any fiber of $\mathcal{X}\to C$
over a point $p\in C$ of good reduction there always exists a very
free curve defined over $\kappa(p)$ and passing through two given
$\kappa(p)$-points. This implies that the comb $C^*$ may be
constructed over $\mathbb{F}$ when $|\mathbb{F}|$ is large enough.

In addition to the teeth given by very free curves in fibers over
some points of $P$, many other teeth have to be added to the comb
$C^*$ such that the condition on the vanishing of higher cohomology
in task (A) holds. Let $m$ be the number of these teeth. The
smallest value for $m$ depends on the normal bundle $N_{\sigma}$ and
the number $\deg_CP$, as we remarked after the statement of
Lemma$\;$\ref{lemma3p3}. By the Lang-Weil estimate (\cite{LW},
Theorem$\;$1), we may assume that $|\mathbb{F}|$ is greater than a
lower bound depending on $\deg\,C\,,\,\deg_CP$ and $|S|$ so that we
can find $m$ $\mathbb{F}$-rational points outside $P\cup S$. We then
construct the $m$ teeth in fibers over these points. The total
number of teeth of the comb $C^*$ is at most $\deg_CP+m$. By
\cite{KoSz}, Theorem$\;$16 and the Lang-Weil estimate (cf.
\cite{KoSz}, p.258, Proof of Theorem 16 implies Theorem 2), we may
choose each tooth with degree bounded from above in terms of $\dim
X\,,\deg\,X$. Thus the degree $d$ of $C^*$ is bounded from above in
terms of $\deg_CP\,,\,m\,,\, \dim X\,,\,\deg\,X$ and
$\deg\,\sigma(C)$.

In the parameter space $M_d$ of degree $d$ curves in $\mathcal{X}$
containing $\{\,r_i\,\}$, the comb $C^*$ corresponds to a smooth
$\mathbb{F}$-rational point $[C^*]$. Choose a smooth geometrically
irreducible curve $T$ in $M_d$ passing through $[C^*]$ such that
deformations given by points in $W=T\setminus\{\,[C^*]\,\}$ are all
smooth curves. Then $W$ is contained in a unique geometrically
irreducible component $W'$ of the subspace $H_d$ parametrizing
smooth degree $d$ curves in $\mathcal{X}$ containing $\{\,r_i\,\}$.
By \cite{KoSz}, Proposition$\;$20, the basic projective invariants
of $W'$ are bounded in terms of $\dim \mathcal{X}$,
$\deg\,\mathcal{X}$, $d$ and $\deg_CP$. In our situation,
$\dim\mathcal{X}$ and $\deg\,\mathcal{X}$ can be expressed in terms
of $\dim X$, $\deg\,C$ and $\deg\,X$. Taking into account all these
things, we finally get, using the Lang--Weil theorem, a lower bound
\[
B=B(N_{\sigma}\,,\, |S|\,,\,
\deg_CP\,,\,\deg\,C\,,\,\deg\,\sigma(C)\,,\,\deg\,X\,,\,\dim\,X)
\]having the property that $W'(\mathbb{F})\neq\emptyset$ whenever $|\mathbb{F}|>B$. Let $\mathcal{C}_w$ be a deformed curve
corresponding to a rational point of $W'$. Deformed curves given by
points of $W$ have the same numerical properties as $C^*$, so does
the curve $\mathcal{C}_w$. In particular, the intersection number of
$\mathcal{C}_w$ with a general fiber of $\mathcal{X}\to C$ is one.
Therefore, the curve $\mathcal{C}_w$ gives a section $s:
C\to\mathcal{X}$ as we want.
\end{proof}

A drawback of the result in Proposition$\;$\ref{prop1p4} is that the
lower bound $B$ is ineffective and depends on many non-intrinsic
objects. However, from the above proof we see that the dependence of
$B$ on the normal bundle $N_{\sigma}$ relies only on its
cohomological property, so  the value of $B$ will not change when we
go over to an extension of the ground field $\mathbb{F}$.

Let us now start the proof of Theorem$\;$\ref{thm1p5}.  Basically,
it amounts to showing that it is really possible to come down to the
ground field $\mathbb{F}$ from a big enough extension $F_n$ over
which the result may be guaranteed by Proposition$\;$\ref{prop1p4}.
Using geometry of cubic hypersurfaces we will do this with
$F_n/\mathbb{F}$ a tower of quadratic extensions. The basic ideas go
back to Swinnerton-Dyer, who obtained similar results for cubic
surfaces over number fields (cf. \cite{SD}, Theorem$\;$5 and its
corollary).

Let $V\subseteq\mathbb{P}^N_k$ be a cubic hypersurface over a field
$k$ and let $k_1/k$ be a quadratic extension. Denote by $\sigma$ the
nontrivial element in the Galois group $\mathrm{Gal}(k_1/k)$. We
define a ``dashed arrow''
\begin{equation}
V(k_1) \dashrightarrow V(k)
\end{equation}as follows: For a point $v_1\in V(k_1)$, the line
$\ell(v_1\,,\,{}^{\sigma}{v}_1)$ joining $v_1$ and its conjugate
point ${}^{\sigma}{v}_1$ will generally intersect $V$ at only one
more point $v_0$, which is $k$-rational. Whenever this is
well-defined, we associate to $v_1\in V(k_1)$ the third intersection
point $v_0\in V(k)$ of the line $\ell(v_1\,,\,{}^{\sigma}{v}_1)$ and
$V$. In what follows, we will always mean the dashed arrow is
defined at $v_1$ and sends it to $v_0$ when we write
``$v_1\dashrightarrow v_0$''.

\begin{lemma}\label{lemma3page5}
Let $A$ be a discrete valuation ring with field of fractions $K$ and
residue field $F$. Let $X_K\subseteq\mathbb{P}^N_K$ be a smooth
cubic hypersurface and $\mathcal{X}_A\subseteq\mathbb{P}^N_A$ its
scheme-theoretic closure in $\mathbb{P}^N_A$. Suppose that the
special fiber $X_F\subseteq\mathbb{P}^N_F$ is a smooth cubic
hypersurface. Let $F_1/F$ be a quadratic extension and $K_1/K$ an
unramified quadratic extension having residue field $F_1/F$. Then
the diagram
\[
\xymatrix{
  X_K(K_1) \ar@{.>}[d]_{} \ar[r]^{} & X_F(F_1) \ar@{.>}[d]^{} \\
  X_K(K) \ar[r]^{} & X_F(F)   }
\]is commutative in the following sense:

If $x_1\in X_K(K_1)$ specializes to $\widetilde{x}_1\in X_F(F_1)$
and $\widetilde{x}_1\dashrightarrow \widetilde{x}_0\in X_F(F)$, then
the dashed arrow $X_K(K_1)\dashrightarrow X_K(K)$ is defined at
$x_1$ and the image $x_0\in X_K(K)$ of $x_1$ specializes to
$\widetilde{x}_0\in X_F(F)$.
\end{lemma}
\begin{proof}
Let $\sigma$ be the nontrivial element of the Galois group
$\mathrm{Gal}(K_1/K)=\mathrm{Gal}(F_1/F)$. The hypothesis
$\widetilde{x}_1\dashrightarrow \widetilde{x}_0$ means that the
$F$-line $L_F=\ell(\widetilde{x}_1\,,\,{}^{\sigma}\widetilde{x}_1)$
joining $\widetilde{x}_1$ and ${}^{\sigma}\widetilde{x}_1$ is not
contained in $X_F$ and the intersection
$\ell(\widetilde{x}_1\,,\,{}^{\sigma}\widetilde{x}_1)\cap X_F$ has
exactly 3 points, the third being $\widetilde{x}_0\in X_F(F)$. Since
$K_1/K$ is unramified, it is easily verified that ${}^{\sigma}x_1$
specializes to ${}^{\sigma}\widetilde{x}_1$. Let
$L_K=\ell(x_1\,,\,{}^{\sigma}x_1)$ be the $K$-line joining $x_1$ and
${}^{\sigma}x_1$. Let $\mathcal{L}_A\subseteq\mathbb{P}^N_A$ be its
scheme-theoretic closure in $\mathbb{P}^N_A$.

{\bf Claim.} $\mathcal{L}_A$ is an $A$-line in $\mathbb{P}^N_A$.

We leave aside the proof of the claim for a moment and continue the
proof of the lemma.

The special fiber $\mathcal{L}_A\times_AF$ is thus a line containing
$\widetilde{x}_1$ and ${}^{\sigma}\widetilde{x}_1$, whence
$\mathcal{L}_A\times_AF=L_F=\ell(\widetilde{x}_1\,,\,{}^{\sigma}\widetilde{x}_1)$.
We have $L_K\nsubseteq X_K$ since $L_F\nsubseteq X_F$. Hence
$L_K\cap X_K=\{\,x_1\,,\,{}^{\sigma}x_1\,,\,x_0\,\}$ with $x_0\in
X_K(K)$, showing that $x_1\dashrightarrow x_0$. To prove $x_0$
specializes to $\widetilde{x}_0$, it suffices to show that the
scheme-theoretic closure of $Z_K:=L_K\cap
X_K\subseteq\mathbb{P}^N_K$ in $\mathbb{P}^N_A$ is equal to
$\mathcal{Z}_A=\mathcal{L}_A\cap
\mathcal{X}_A\subseteq\mathbb{P}^N_A$. By \cite{Gr3}, (IV.2.8.5), we
need only show that $\mathcal{Z}_A$ is $A$-flat. Note that
$\mathcal{Z}_A\to\mathrm{Spec}\, A$ is proper and both the generic
fiber $Z_K/K$ and the special fiber $Z_F=L_F\cap X_F/F$ are finite.
By Chevalley's theorem, $\mathcal{Z}_A$ is finite over $A$. Since
the generic fiber $Z_K/K$ and the special fiber $Z_F/F$ have the
same length, $\mathcal{Z}_A$ is flat over $A$ by \cite{Har}, p.174,
Lemma$\;$8.9.

We finish by giving the proof of our claim. This is essentially an
easy consequence of \cite{Gr3}, (IV.2.8.1.1).

The embedding $L_K\subseteq\mathbb{P}^N_K$ corresponds to a
surjective $K$-homomorphism $\varphi: K^{N+1}\to K^2$. Let $M$ be
the image of the composite map $A^{N+1}\hookrightarrow
K^{N+1}\overset{\varphi}{\longrightarrow}K^2$ and let $\varphi_A:
A^{N+1}\to M$ be the induced map. Then the surjection $\varphi_A:
A^{N+1}\to M$ corresponds to the embedding
$\mathcal{L}_A\subseteq\mathbb{P}^N_A$. We need prove that $M$ is
free of rank 2 over $A$. In fact, there is an induced commutative
diagram with exact rows
\[
\xymatrix{
  A^{N+1} \ar[d]_{} \ar[r]^{\varphi_A} & M \ar[d]_{} \ar[r]^{} & 0 \\
  K^{N+1} \ar[r]^{\varphi} & K^2 \ar[r]^{} & 0   }
\]As an $A$-submodule of $K^2$, $M$ is torsion-free and hence free
over the discrete valuation ring $A$. The natural homomorphism
$M\otimes_AK\to K^2$ is bijective, whence $\mathrm{rank}_AM=2$. This
completes the proof.
\end{proof}

We shall now make use of the following lemma due to Koll\'ar.

\begin{lemma}[\cite{Kol4}, Lemma 9.4]\label{lemmaAdded1}
Let $\mathbb{F}_q$ be a finite field of cardinality $q\ge
11$.\footnote{In \cite{Kol4}, the proof of Lemma$\;$9.4 holds only
for $q\ge 11$.} Then for any smooth cubic hypersurface
$X\subseteq\mathbb{P}^N$ defined over ${\mathbb{F}_q}$ with $N\ge
2$, the dashed arrow
\[
X(\mathbb{F}_{q^2})\dashrightarrow X(\mathbb{F}_q)
\]is surjective, i.e., for every $p\in X(\mathbb{F}_q)$, there is a point $x\in
X(\mathbb{F}_{q^2})\setminus X(\mathbb{F}_q)$ such that the line
$\ell(x\,,\,{}^{\sigma}{x})$ joining $x$ and its conjugate point
${}^{\sigma}{x}$ satisfies
\[
\ell(x\,,\,{}^{\sigma}{x})\cap
X=\{\,p\,,\,x\,,\,{}^{\sigma}{x}\,\}\,.
\]
\end{lemma}

\begin{proof}[Proof of Theorem$\;\ref{thm1p5}$]We consider the case
where $P=\{\,p\,\}$ consists of a single point. The general case may
be treated in the same way without essential difference.

Let $\kappa=\kappa(p)$ be the residue field of the point. We fix an
algebraic closure $\overline{\mathbb{F}}$ of $\mathbb{F}$ and let
$F_n\,,\,\kappa_n$, $n\ge 0$ be the subfields of
$\overline{\mathbb{F}}$ determined by the following conditions:
\[
F_0=\mathbb{F}\,,\;\; \kappa_0=\kappa\,,\; F_n\subseteq
\kappa_n\,,\;\text{ and }\; [F_{n+1}: F_n]=[\kappa_{n+1}:
\kappa_n]=2\,,
\]for all $n\ge 0$. Let $K_n=F_n(C_{F_n})$ be the function field of
the curve $C_{F_n}=C\times_{\mathbb{F}}F_n$ defined over $F_n$. Then
for each $n$ the following diagram is commutative in the sense
described as in Lemma$\;$\ref{lemma3page5}:
\[
\xymatrix{
  X(K_{n+1}) \ar@{.>}[d]_{} \ar[r]^{} & \mathcal{X}_p(\kappa_{n+1}) \ar@{.>}[d]^{} \\
  X(K_n) \ar[r]^{} & \mathcal{X}_p(\kappa_n)   }
\]Suppose that
$|\kappa(p)|>\phi(N)$. Then, by Lemma$\;$\ref{lemmaAdded1}, starting
from any point $r_0\in \mathcal{X}_p(\kappa_0)$ we can find
successively a sequence of points $r_i\in \mathcal{X}_p(\kappa_i)$
such that for each $i$, $r_{i+1}\dashrightarrow r_i$ via the dashed
arrow $\mathcal{X}_p(\kappa_{i+1})\dashrightarrow
\mathcal{X}_p(\kappa_i)$. Due to Proposition$\;$\ref{prop1p4}, we
may choose $n$ big enough so that the specialization map
$X(K_{n+1})\to \mathcal{X}_p(\kappa_{n+1})$ is surjective. Pick a
point $s_{n+1}\in X(K_{n+1})$ that specializes to
$r_{n+1}\in\mathcal{X}_p(\kappa_{n+1})$. By
Lemma$\;$\ref{lemma3page5}, we obtain points $s_i\in X(K_i)$ for
$i=n+1,\,n\,\dotsc, 0$ such that
\[
s_{n+1}\dashrightarrow s_n \dashrightarrow\cdots \cdots
\dashrightarrow s_1\dashrightarrow s_0
\]and each $s_i$ specializes to $r_i$. In particular, there exists a
point $s_0\in X(K)$ which specializes to
$r_0\in\mathcal{X}_p(\kappa(p))$.
\end{proof}

\subsection{Weak approximation at one place}

Our proof of Theorem$\;$\ref{thm1p6} follows the method of
Swinnerton-Dyer. In particular, the following result due to him will
be needed.

\begin{thm}[\cite{SD}, Theorem$\;$4]\label{thm4ofSD}
Let $K$ be a global field, $K_v$ the completion of $K$ at a
nonarchimedean place $v$ and $k$ the residue field of $v$. Let
$V\subseteq\mathbb{P}^3_K$ be a smooth cubic surface whose reduction
$\widetilde V\subseteq \mathbb{P}^3_k$ at $v$ is also smooth. Let
$\widetilde P\in \widetilde V(k)$ be the reduction of a point $P\in
V(K)$.

Suppose that there is a point $R\in V(K)$ whose reduction
$\widetilde{R}\in\widetilde V(k)$ has the following properties:

$(1)$ the line $\widetilde{P}\widetilde R$ intersects $\widetilde V$
at exactly three distinct points;

$(2)$ no geometric line on $\widetilde V$ passes through $\widetilde
R$;

$(3)$ letting $\gamma=T_{\widetilde R}\widetilde V\cap\widetilde V$,
there exist two distinct points $t_1,\,t_2\in\gamma(k)$ such that
for each $i=1,\,2$, $t_i\neq\widetilde R$ and the line $\widetilde P
t_i$ intersects $\widetilde V$ at three distinct points and
\[
T_{t_1}\gamma\cap T_{t_2}\gamma\cap T_{\widetilde P}\widetilde
V=\emptyset\,.
\]

Then there is an $R$-equivalence class $\mathscr{C}$ in $V(K)$ such
that for every point $Q^*\in V(K_v)$ that specializes to $\widetilde
P\in \widetilde V(k)$ and every $v$-adic open neighborhood
$\mathscr{U}_v$ of $Q^*$, one has $\mathscr{C}\cap
\mathscr{U}_v\neq\emptyset$.
\end{thm}

Swinnerton-Dyer stated the above theorem for a number field $K$. But
his proof works for function fields as well. (In his paper
\cite{SD}, Swinnerton-Dyer deduced this theorem from his Lemma$\;$8
and Theorem$\;$3. In proving these two results the key ingredients
have been nonarchimedean analysis plus some geometric arguments
involving tangent spaces and intersection theory, which do not
actually depend on special properties of number fields as one may
verify with more or less patience.)

\

In the proof of the next result, some geometric arguments have been
used in \cite{SD} and \cite{Kol4}. The idea being very similar, here
we still include a  proof because of a slightly greater generality
and some useful explicit estimates.

\begin{prop}\label{prop5page9}
There is a function $\phi: \mathbb{N}\to \mathbb{N}$ having the
following property:

Let $\mathbb{F}_q$ be a finite field, $N\ge 2$, and
$X\subseteq\mathbb{P}^N$ a smooth cubic hypersurface over
$\mathbb{F}_q$. If $q=|\mathbb{F}_q|>\phi(N)$, then for every $p\in
X(\mathbb{F}_q)$, there is a point $x\in X(\mathbb{F}_q)$ such that
the line $\ell(x,\,p)$ joining $x$ and $p$ intersects $X$ at exactly
three distinct $\mathbb{F}_q$-points.
\end{prop}
\begin{proof}
Fix $p\in X(\mathbb{F}_q)$. Set
\[
B=\{\,x\in X(\mathbb{F}_q)\,|\,x\neq p\,\text{ and }\;
\ell(x,\,p)\nsubseteq T_pX\,\}\,,
\]where $T_pX\subseteq\mathbb{P}^N$ denotes the tangent hyperplane to $X$ at
$p$. For a line $\ell\subseteq\mathbb{P}^N$ defined over
$\mathbb{F}_q$ not contained in $X$ and a point $x\in (\ell\cap
X)(\mathbb{F}_q)$, we denote by $(\ell\cdot X)_x$ the intersection
multiplicity of $\ell$ with $X$ at the point $x$. Then $B=B_1\cup
B_2$, where for $i=1\,,2\,$, the subset $B_i$ consists of points
$x\in B$
 such that $(\ell(x,\,p)\cdot X)_x=i$. Let $b=\#B$
and $b_i=\#B_i\,,i=1\,,2$. We have
\[
\begin{split}
b=&b_1+b_2=\#(X(\mathbb{F}_q)\setminus\{\,p\,\})-\\
&\#\{\,\mathbb{F}_q\text{-line }\,\ell\subseteq T_pX\,|\,p\in\ell\;
\text{ and }\;(\ell\cdot X)_p=2\,\}\,
\end{split}\]whence
\begin{equation}\label{eq3page10}
\#X(\mathbb{F}_q)-1-\frac{q^{N-1}-1}{q-1}\le b=b_1+b_2\le
\#X(\mathbb{F}_q)-1\,.
\end{equation}We want to show that for $q$ large enough,
\[
\#X(\mathbb{F}_q)-1-\frac{q^{N-1}-1}{q-1}>b_2\,
\]so that $B_1\neq\emptyset$. Observe that
\[
B_2=\{\,x\in X(\mathbb{F}_q)\,|\, x\notin T_pX\,,\;\; p\in
T_xX\,\}\,.
\]Suppose that $X$ is defined by a cubic form $\varphi\in \mathbb{F}_q[T_0,\,\cdots, T_N]$
and put $\varphi'_0=\frac{\partial \varphi}{\partial T_0}$. We may
assume that $p=(1:0:\cdots:0)$. Letting $Y\subseteq\mathbb{P}^N$
denote the subscheme defined by the equations
$\varphi=\varphi'_0=0$, we have clearly $B_2\subseteq
Y(\mathbb{F}_q)$, whence $b_2\le \# Y(\mathbb{F}_q)$.

By the Lang-Weil theorem,
\[
\#X(\mathbb{F}_q)=q^{N-1}+O(q^{N-3/2})
\]where the term $O(q^{N-3/2})$ depends only on $N$ and $q$. Using \cite{Del}, Th\'{e}or\`{e}me 8.1,
we can get an even more explicit estimate:
\[
\left|\#X(\mathbb{F}_q)-\frac{q^{N}-1}{q-1}\right|\le \beta_N\cdot
q^{\frac{N-1}{2}}
\]where $\beta_N$ denotes the $(N-1)$-th Betti number of a smooth cubic
hypersurface in $\mathbb{P}^N_{\mathbb{C}}$. In view of
\eqref{eq3page10}, it is now sufficient to show that
\begin{equation}\label{eq4page11}
\# Y(\mathbb{F}_q)=O(q^{N-2})
\end{equation}with the term $O(q^{N-2})$ depending only on $N$ and
$q$. An easy verification shows that the form $\varphi'_0$ cannot be
identically zero on $X$ since $X$ is smooth and geometrically
irreducible. It follows that $\dim Y<\dim X=N-1$ and thus
\eqref{eq4page11} follows from the classical Lang-Weil bound. In our
case, it is also possible to give an explicit bound for
$\#Y(\mathbb{F}_q)$. For instance, \cite{Schm}, Chapter 4,
Lemma$\;$3.3 gives
\begin{equation}\label{eq5p4}
\# Y(\mathbb{F}_q)\le \frac{12q^{N-1}}{q-1}\,.
\end{equation}The proposition is
thus proved.\qed
\end{proof}

\begin{remark}\label{remark5p4}
According to the above proof (cf. \eqref{eq3page10} and
\eqref{eq5p4}), for a function $\phi: \mathbb{N}\to \mathbb{N}$ to
have the property in Proposition$\;$\ref{prop5page9}, one sufficient
condition is
\[
\frac{12q^{N-1}}{q-1}<
q^{N-1}-1-\beta_Nq^{\frac{N-1}{2}}\,,\;\;\;\forall\; N\ge
2\,,\;\forall\; q>\phi(N)\,.
\]So one can compute explicitly at least one possible value of $\phi(N)$ for each $N$.
For instance, using $\beta_3=7$ one verifies easily that one may
take $\phi(3)=20$.
\end{remark}

\begin{remark}\label{remarkAdded2}
Letting $\mathscr{L}_N$ be the set of $\mathbb{F}_q$-lines
$\ell\subseteq\mathbb{P}^N$ such that
\[
\ell\cap X=\{p\,,\,x\,,\,{}^{\sigma}x\}\,,\quad \text{ with }\;x\in
X(\mathbb{F}_{q^2})\setminus X(\mathbb{F}_q)\,,
\]our proof of Proposition \ref{prop5page9} can also give a lower
bound depending on $N$ (and $q$) for the cardinality
$\#\mathscr{L}_N$. But as a sufficient condition for
$\mathscr{L}_N\neq\emptyset$, Koll\'ar's lemma (Lemma
\ref{lemmaAdded1}) is obviously better.
\end{remark}

\begin{lemma}\label{lemma5p6}
Let $\widetilde V$ be a smooth cubic surface over a finite field $k$
of characteristic not dividing $6$. If $q=|k|>47$, then for any
$\widetilde P\in \widetilde V(k)$, one can find a point $\widetilde
R\in\widetilde V(k)$ satisfying all the conditions $(1),\,(2)$ and
$(3)$ of Theorem$\;\ref{thm4ofSD}$.
\end{lemma}

\begin{proof}
Swinnerton-Dyer has noticed this fact as he remarked in \cite{SD},
p.379, lines 7--10. However, an explicit proof seems not included
there. So we give our own proof here.

Let $B_1$ be the set of points $\widetilde R\in\widetilde V(k)$ for
which the condition (1) holds. We know from the proof of
Proposition$\;$\ref{prop5page9} (cf. \eqref{eq3page10} and
\eqref{eq5p4}) that the cardinality $b_1=\#B_1$ satisfies
\[
b_1\ge
\frac{q^3-1}{q-1}-\beta_3q-\left(1+\frac{13q^2-1}{q-1}\right)=q^2-6q-\frac{13q^2-1}{q-1}\,.
\]
Since the union of all geometric lines on $\widetilde V$ contains at
most $27(q+1)$ $k$-points, the subset $B'_1\subseteq B_1$ consisting
of points $\widetilde{R}\in\widetilde V(k)$ which have both the
properties (1) and (2) has cardinality
\[
b'_1=\#B_1'\ge b_1-27(q+1)\ge q^2-46q-52\,.
\]
For $q>47$, we can always find a point $\widetilde R\in\widetilde
V(k)$ having the properties (1) and (2). Then the intersection
$\gamma=T_{\widetilde R}\widetilde V\cap \widetilde V$ is a
geometrically irreducible plane cubic curve. Such a curve has at
least $q-2$ smooth rational points. The points $t\in \widetilde
V(k)$ such that $\widetilde P\in T_t\widetilde V$ all lie on a
quadratic $Y$ and the intersection $\gamma\cap Y$ hat at most 6
points. Thus, when $q>47$, we can always find points $t_1\neq
t_2\in\gamma(k)$ with $t_i\neq \widetilde R$ such that the line
$\widetilde Pt_i$ intersects $\widetilde V$ at 3 distinct points for
each $i=1,\,2$. It remains to show that $t_1,\,t_2$ may be so chosen
that the intersection
\[
T_{t_1}\gamma\cap T_{t_2}\gamma\cap T_{\widetilde P}\widetilde
V=T_{t_1}\widetilde V\cap T_{t_2}\widetilde V\cap
\left(T_{\widetilde R}\widetilde V\cap T_{\widetilde P}\widetilde
V\right)
\]is empty. This is a consequence of the Lemma$\;$\ref{lemma5p7}
below.
\end{proof}

\begin{lemma}\label{lemma5p7}
Let $F$ be any field of characteristic not dividing $6$,
$\gamma\subseteq\mathbb{P}^2_F$ a geometrically irreducible plane
cubic curve having a nonsmooth point $R\in\gamma(F)$. Then for any
$S\in\mathbb{P}^2(F)$, the set
\[
\mathscr{T}:=\{\,t\in\gamma(F)\,|\, \text{ the line $St$ joining $S$
and $t$ is tangent to $\gamma$ at }\;t \,\}
\]has cardinality $\le 6$.
\end{lemma}
\begin{proof}
Take coordinates such that $R=(1:0:0)$. Since $\gamma$ passes
through $R$ and is not smooth there, the equation of $\gamma$ has
the following form
\[
\varphi=T_0q(T_1\,,\,T_2)+c(T_1\,,\,T_2)
\]where $q$ is a quadratic form and $c$ is a cubic form. We may
assume $S\neq R$ and take coordinates such that $S=(0:1:0)$. Then
$\mathscr{T}$ is the set of rational points of the subscheme
$Z\subseteq\mathbb{P}^2_F$ defined by
\[
\varphi=\frac{\partial \varphi}{\partial T_1}=0\,.
\]Using the hypothesis on the characteristic of $F$, one
concludes by an easy computation that $\frac{\partial
\varphi}{\partial T_1}$ is not identically zero because of the
geometrical irreducibility of $\gamma$. It then follows that
$\#\mathscr{T}=\# Z(F)\le 6$.
\end{proof}

\begin{proof}[Proof of Theorem$\;\ref{thm1p6}$]
We apply Theorem$\;$\ref{thm4ofSD} with $V$ equal to the generic
fiber $X$ and $\widetilde V$ equal to the special fiber
$\mathcal{X}_p$. By Theorem$\;$\ref{thm1p5}, the hypothesis implies
that lifting a point in $\widetilde V(k)$ to a point in $V(K)$ is
always possible. The result then follows easily from
Theorem$\;$\ref{thm4ofSD} and Lemma$\;$\ref{lemma5p6}.
\end{proof}

\noindent \emph{Acknowledgements.} The author wishes to thank the
referee for valuable comments, especially for pointing out the
possibility of generalizing the first version of
Theorem$\;$\ref{thm1p5} using the reference \cite{Kol4}. The author
also thanks Prof.\! Jean-Louis Colliot-Th\'{e}l\`{e}ne for many
valuable discussions and comments. Some of the problems solved here
arose out of conversations which Jean-Louis Colliot-Th\'{e}l\`{e}ne
and Brendan Hassett had in 2005, and the approach we follow is also
inspired by their discussions. Thanks also go to Prof.\! J\'anos
Koll\'ar for useful discussions about cubic hypersurfaces over
finite fields.

\addcontentsline{toc}{section}{

\end{document}